\documentclass{mystyle} 
\usepackage{latexsym} 
\usepackage{alltt}
\usepackage[all]{xy} 
\CompileMatrices 
\xyoption{frame}
\usepackage{amssymb}
\usepackage{amscd} 
\usepackage{amsmath}
\usepackage{theorem} 
\theoremstyle{break} \newtheorem{theorem}{Theorem}[section]
\theoremstyle{break} \newtheorem{corollary}[theorem]{Corollary}
\theoremstyle{break} \newtheorem{lemma}[theorem]{Lemma}
\theoremstyle{break} \newtheorem{proposition}[theorem]{Proposition}
\theoremstyle{plain} {\theorembodyfont{\rmfamily} \newtheorem{example}[theorem]{Example}}
\theoremstyle{plain} \newtheorem{definition}[theorem]{Definition}

\begin{document} 
\begin{frontmatter} 
\title{Some properties of finite meadows}
\author[label1]{Inge Bethke\thanksref{email1}\corauthref{cor}}
\author[label1]{Piet Rodenburg\thanksref{email2}}
\corauth[cor]{Corresponding author. Address: Kruislaan 403, 1098 SJ Amsterdam, The Netherlands }
\thanks[email1]{E-mail:  \texttt{inge@science.uva.nl}}
\thanks[email2]{E-mail:  \texttt{pietr@science.uva.nl}}
\address[label1]{University of Amsterdam, Faculty of Science, Section Theoretical
Software
Engineering (former Programming Research Group)} 

\begin{keyword}
combinatorial problems, data structures, specification languages.
\end{keyword}

\end{frontmatter}
\section{Introduction}
In abstract algebra, a field is a structure with total operations of addition, subtraction
and multiplication. Moreover, every element has a multiplicative inverse|except 0. In a field, the rules 
hold
which are familiar from the arithmetic of ordinary numbers.
The prototypical example is the field of rational numbers. Fields  can be specified by the axioms for 
commutative rings with identity element
($\mathit{CR}$, see Table \ref{meadowaxioms}), and the negative conditional formula 
\[
x\neq 0 \rightarrow x\cdot x^{-1}=1,
\]
which is difficult to apply and automate in formal reasoning.

The theory of fields is a very active area which is not only of great theoretical interest but has also found
applications both within mathematics|combina\-to\-rics and algorithm analysis|as well as in engineering sciences 
and, in particular, in coding theory and 
sequence design. Unfortunately, since fields are not 
axiomatized by equations only, Birkhoff's Theorem fails, i.e. fields do not
constitute a variety: they are not closed under products, subalgebras, and homomorphic images.

In \cite{BT07}, it is proved that there exists a finite equational specification under initial algebra 
semantics|without hidden functions|of the rational numbers with field operations that are all total.
Subsequent investigations led to the concept of \emph{meadows} which are very similar to fields|the 
considerable difference being that meadows do form a variety.

A meadow is a commutative ring with identity element ($\mathit{CR}$) equipped with a total unary operation $\ ^{-1}$, 
\emph{inversion},
which satisfies the equation for \emph{reflection ($\mathit{Ref}$)} and the \emph{restricted
inverse law ($\mathit{Ril}$)}. That is, a meadow is specified by the set of axioms in Table \ref{axioms}.
\begin{table}\label{axioms}
\[
\begin{array}{lrcl}
\hline
(\mathit{CR})\hspace{1 cm}&(x+y)+z &=&x+(y+z)\\
&x+y&=&y+x\\
&x +0 &=& x\\
&x+(-x)&=&0\\
&(x \cdot y)\cdot z&=& x\cdot (y \cdot z)\\
&x\cdot y &=&y \cdot x\\
&x \cdot 1 &=& x\\
&x\cdot (y + z) &=& x\cdot y + x \cdot z\\
(\mathit{Ref}) &(x^{-1})^{-1} &=& x\\
(\mathit{Ril}) & x \cdot (x \cdot x^{-1})& = &x\\
\hline
\end{array}
\]
\caption{Specification of the theory of meadows \label{meadowaxioms}}
\end{table}
All fields and products of fields can be viewed as meadows|basically by stipulating $0^{-1}=0$|but not conversely.
 Also, every commutative Von Neumann regular ring (see e.g.\ \cite{G79}) can be expanded
to a meadow (cf.\ \cite{BHT07}).

\begin{example}\label{example}
Consider the ring $\mathbb{Z}/10\mathbb{Z}$ with elements $\{0,1,2, \ldots , 9\}$ where arithmetic is performed modulo $10$. Here
\[
\begin{array}{rclccrcl}
(0)^{-1} &=& 0&\hspace{1 cm}&
(1)^{-1} &=& 1\\
(2)^{-1} &=& 8&&
(3)^{-1} &=& 7\\
(4)^{-1} &=& 4&&
(5)^{-1} &=& 5\\
(6)^{-1} &=& 6&&
(7)^{-1} &=& 3\\
(8)^{-1} &=& 2&&
(9)^{-1} &=& 9\\
\end{array}
\]
Since the inversion is an involution which also satisfies $\mathit{Ril}$, this ring is also a meadow.
\end{example}
The aim of this note is to describe the structure of finite meadows. 
We will show that the class of finite meadows is the closure
of the class of finite fields under finite products. As a corollary, we obtain a unique representation of
minimal meadows in terms of prime fields.

\section{Decomposition of finite meadows}
In \cite{B44} it is proved that every commutative regular ring in the sense of von Neumann
is a subdirect union of fields.
In this section we show that every finite meadow is a direct product of 
finite fields. Part of the proof is also known from the theory of rings:
under certain conditions|also met in our case|a ring $R$ can be decomposed as 
$R=e_1\cdot R\cdot e_1 \oplus \ldots \oplus e_n\cdot R\cdot e_n$ where $\{e_1, \ldots ,
e_n\}$ is the set of mutually orthogonal minimal idempotents in $R$ (see e.g.\ \cite{D06}).
\begin{definition}
Let $M$ be a meadow.
\begin{enumerate}
\item An element $e\neq 0$ in $M$ is an \emph{idempotent} if $e \cdot e = e$. 
\item If $e, e'\in M$ are idempotents then we write $e \leq e'$ if $e\cdot e'=e$.
\item An idempotent $e\in M$ is \emph{minimal} if for every idempotent
$e'\in M$, 
\[e'\leq e \Rightarrow e'=e.\] 
\end{enumerate}
\end{definition}
\begin{proposition}\label{meadow}
Let $M$ be a meadow and $e\in M$ an idempotent.
Then
\begin{enumerate}
\item $e=e^{-1}$
\item $e\cdot M$ is a meadow with multiplicative identity element $e$.
\item If $e$ is minimal then $e\cdot M$ is a field with multiplicative identity element $e$.
\end{enumerate}
\end{proposition}
{\bf Proof:} 
\begin{enumerate}
\item Since in every meadow inversion distributes over multiplication (see \cite{BHT07}) we have
\[
e = e\cdot e\cdot e^{-1}=e\cdot e^{-1}= e \cdot (e\cdot e)^{-1}= 
e\cdot e^{-1}\cdot e^{-1}= e^{-1}.
\]
\item Since idempotents are self-inverse $e\cdot M$ is closed under $+, \cdot, \-^{-1}$ and clearly satisfies the
 axioms for meadows.
\item Since $e\cdot M$ is a meadow with multiplicative identity element $e$, it suffices to prove that 
$(e\cdot m)\cdot (e\cdot m)^{-1} = e$ for every $e\cdot m \neq 0$. Thus let $e\cdot m$ be a nonzero element. 
Then $(e\cdot m)\cdot (e\cdot m)^{-1} \neq 0$ because otherwise
 \[
e\cdot m = (e\cdot
 m)\cdot (e \cdot m) \cdot (e\cdot m)^{-1} = 0.
 \]
 Moreover,
\[
(e\cdot m)\cdot (e\cdot m)^{-1}\cdot (e\cdot m)\cdot (e\cdot m)^{-1} = (e\cdot m)\cdot (e\cdot m)^{-1}.
\]
So $(e\cdot m)\cdot (e\cdot m)^{-1}$ is an idempotent. Hence, since
\[
e\cdot (e\cdot m)\cdot (e\cdot m)^{-1} = (e\cdot m)\cdot (e\cdot m)^{-1}
\]
and $e$ is minimal we have $(e\cdot m)\cdot (e\cdot m)^{-1} = e$.
\end{enumerate}
$\qed$

The main properties of idempotents are summarized in the following proposition.
\begin{proposition}\label{properties}
Let $M$ be a meadow. 
\begin{enumerate}
\item $\leq$ is a partial order on the idempotents.
\item If $e,e'\in M$ are idempotents and $e\cdot e'\neq 0$ then $e\cdot e'$ is also an idempotent.
\item If $e,e'\in M$ are  idempotents and $e < e'$ then $e'-e$ is also an idempotent.
\end{enumerate}
\end{proposition}
{\bf Proof}:
\begin{enumerate}
\item Clearly $\leq$ is reflexive. If $e\leq e'$ and $e'\leq e''$ then 
\[
e\cdot e'' = (e\cdot e') \cdot e'' = e\cdot (e' \cdot e'') = e\cdot e' =e.
\]
Therefore the relation is transitive. Finally, if $e\leq e'$ and $e'\leq e$ then
\[
e = e\cdot e' = e'\cdot e=e.'
\]
Thus $\leq$ is also antisymmetric.
\item We multiply $e\cdot e'$ with itself:
$
(e\cdot e')\cdot (e\cdot e') = (e\cdot e)\cdot (e'\cdot e') =e\cdot e'.
$
\item We multiply $e'-e$ with itself:
\[
(e'-e)\cdot (e'-e)= e'\cdot e' -e\cdot e' -e'\cdot e +e\cdot e = e' -e -e +e = e'-e.
\]
\end{enumerate}
$\qed$
\begin{definition}
Let $M$ be a meadow and 
$e,e' \in M$. We call $e$ and $e'$ \emph{orthogonal} if $e\cdot e'=0$.
\end{definition}
\begin{proposition}\label{orthogonal}
Let $M$ be a meadow. 
\begin{enumerate}
\item If $e,e'\in M$ are different minimal idempotents then $e$ and $e'$ are orthogonal.
\item If $e,e'\in M$ are orthogonal idempotents then $e+e'$ is an idempotent.
\end{enumerate}
\end{proposition}
{\bf Proof}:
\begin{enumerate}
\item Suppose $e\cdot e'\neq 0$. Then $e\cdot e'$ is an idempotent by Proposition \ref{properties}(2). Moreover,
$
e\cdot e' = e\cdot e \cdot e' = e\cdot e' \cdot e,
$
i.e.\ $e\cdot e' \leq	 e$. Thus $e\cdot e' = e$, since $e$ is minimal. Likewise $e\cdot e' = e'$ and hence
$e=e'$. Contradiction.
\item We multiply again:
\[
(e + e')\cdot (e+ e') = e\cdot e + e \cdot e' + e' \cdot e + e' \cdot e' = e + 0 
+ 0 + e' =e + e'.
\]
Moreover, $(e + e') \cdot e = e\cdot e + e\cdot e'= e$. Hence $e+e'\neq 0$.
\end{enumerate}
$\qed$

We now show that every finite meadow is the direct product of the fields generated by its minimal 
idempotents.
\begin{lemma}\label{som}
Let $M$ be a finite meadow and $\{e_1, \ldots ,e_n\} \subseteq M$ be the set of minimal
idempotents. Then $e_1 + \cdots  + e_n = 1$.
\end{lemma}
{\bf Proof}: Since minimal idempotents are orthogonal we have $e_i\cdot e_j = 0$ 
for $i\neq j$ by Proposition \ref{orthogonal} (1). Therefore for every $1 \leq i < n$, $e_1 + \cdots + e_i$ is an 
idempotent orthogonal with $e_{i+1}$, and hence $e_1 + \cdots  + e_n$ is an idempotent by
Proposition \ref{orthogonal} (2). And therefore $1 - (e_1 + \cdots  + e_n)$ is an idempotent by Proposition
\ref{properties} (3) unless it is 0. Suppose $1 - (e_1 + \cdots  + e_n)$ is an idempotent. Then, since $\leq$ is
a partial order there must be some minimal idempotent $e_i\leq 1 - (e_1 + \cdots  + e_n)$. But
\[
\begin{array}{rcl}
e_i \cdot (1 - (e_1 + \cdots  + e_n))& =& e_i -(e_i\cdot e_1 + \cdots + e_i\cdot e_i + \cdots + e_i\cdot e_n)\\
&= &e_i - (0 + \cdots + e_i\cdot e_i + \cdots + 0)\\
& = &0
\end{array}
\]
Contradiction. Hence $1 - (e_1 + \cdots  + e_n)$ is not an idempotent, i.e. 
\[
1 - (e_1 + \cdots  + e_n)=0
\]
whence $e_1 + \cdots  + e_n = 1$. $\qed$

\begin{theorem}
Let $M$ be a finite meadow and $\{e_1, \ldots ,e_n\} \subseteq M$ the set of minimal
idempotents. Then 
\[
M \cong e_1\cdot M \times \cdots \times e_n \cdot M
\]
\end{theorem}
{\bf Proof}: Because the theory of meadows is equational, we know from universal algebra that 
a direct product of meadows is a meadow.
Thus $e_1\cdot M \times \cdots \times e_n \cdot M$ is a meadow with multiplicative 
identity element
$(e_1, \ldots , e_n)$ and the operations defined componentwise.
Define $h: M \rightarrow e_1\cdot M \times \cdots \times e_n \cdot M$ by
\[
h(m) = (e_1\cdot m, \ldots , e_n\cdot m).
\]
 Then $h$ is a homomorphism. 
Suppose $h(m)=h(m')$. Then for every $1\leq i\leq n$, $e_i\cdot m =e_i\cdot m'$.
Thus
\[
\begin{array}{rcl}
m=1\cdot m&=& (e_1 + \cdots +e_n)\cdot m\\
&=& e_1\cdot m +\cdots +e_n \cdot m \\ 
&=&e_1\cdot m' +\cdots +e_n \cdot m' \\
&= &(e_1 + \cdots +e_n)\cdot m' = 1\cdot m'=m'.
\end{array}
\]
Hence $h$ is injective.
Now let $(e_1\cdot m_1, \ldots , e_n\cdot m_n)\in e_1\cdot M \times \cdots \times e_n \cdot M$ and 
consider $m= e_1\cdot m_1 + \ldots + e_n\cdot m_n$. Then, 
since $e_i$ and $e_j$ are orthogonal for $i\neq j$,
$e_i\cdot m = e_i \cdot m_i$. Thus $h(m)=(e_1\cdot m_1, \ldots , e_n\cdot m_n)$. Whence 
$h$ is also surjective. $\qed$

The order, or number of elements, of finite fields
is of the form $p^n$, where $p$ is a prime number.  Since any two finite 
fields with the same number of elements are isomorphic, there is a naming 
scheme of finite fields that specifies only the order of the field. 
One notation for a finite field|or more precisely, its zero-totalized
expansion, in which inverse is a total operation with $0^{-1}=0$|with $p^n$ elements is $GF(p^n)$, 
where the letters $GF$ stand for \emph{Galois field}. From the above theorem 
it now follows immediately that the class of finite meadows is 
the closure of  the class of Galois fields under finite products.
\begin{corollary}\label{decomp}
Let $\Sigma =\{0,1, +, -, \cdot, \ ^{-1}\}$ and $M$ be a finite $\Sigma$-structure with cardinality $n$. Then $M$ is a meadow
if and only if there are|not necessarily distinct|primes $p_1, \ldots
, p_k$ and natural numbers $n_1, \ldots , n_k$ such that
\[
M \cong GF(p_1^{n_1})\times \cdots \times GF(p_k^{n_k})
\]
and $n= p_1^{n_1} \cdots p_k^{n_k}$.
\end{corollary}
Observe that|as a consequence|meadows of the same size are not necessarily isomorphic:
$GF(4)$ and $GF(2)\times GF(2)$ are both meadows but $GF(4)\not \cong GF(2)\times GF(2)$. However, minimal finite meadows|i.e.\ meadows containing no proper submeadows|of the same size are isomorphic.
\begin{corollary}
\begin{enumerate}
\item Let $M$ be a finite minimal meadow with cardinality $n$. Then there are distinct primes $p_1, \ldots
, p_k$ such that
\[
M \cong \mathbb{Z}/p_1 \mathbb{Z} \times \cdots \times \mathbb{Z}/p_k \mathbb{Z} 
\]
and $n= p_1 \cdots p_k$.
\item Finite minimal meadows of the same size are isomorphic.
\end{enumerate}
\end{corollary}
{\bf Proof}: (2) follows from (1) and
(1) follows from the preceding corollary and the fact that every minimal meadow has $n$ elements, where $n$ is squarefree, 
i.e. its prime factor decomposition is the product of distinct primes (cf. \cite{BHT07}). $\qed$

As an application of Corollary \ref{decomp}, we determine the number of self-inverse
and invertible elements in finite meadows.
\begin{definition}
Let $M=\langle M, 0,1,+,-, \cdot, \ ^{-1}\rangle$ be a meadow and
$m\in M$. Then
\begin{enumerate}
\item $m$ is \emph{self-inverse} if $m = m^{-1}$,
\item $m$ is \emph{invertible} if $m \cdot m^{-1} = 1$,
\end{enumerate}
\end{definition}
So, e.g.\ in $\mathbb{Z}/10\mathbb{Z}$ (see Example \ref{example}) $0,1,4,5,6$ are self-inverse elements,
$1, 3,7,$ are invertibles, and $9$ is both self-inverse and invertible.
\begin{proposition}
Let $M\cong GF(p_1^{k_1})\times \cdots \times GF(p_n^{k_n})$. Then
$M$ has
\begin{enumerate}
\item $2^l\cdot 3^{n-l}$ self-inverses, where $l=\mid\{i\mid 1\leq i\leq n \ \&\  
p_i = 2\ \& \ k_i=1\}\mid$, and 
\item $(p_1^{k_1}-1) \cdots  (p_n^{k_n}-1)$ invertibles.
\end{enumerate}
\end{proposition}
{\bf Proof}: First observe that the number of self-inverses [invertibles] of $M$ is the product of
the number of self-inverses [invertibles] in the Galois fields.\\
(1) Now $m$ is self-inverse in a meadow iff $m^3=m\cdot m\cdot m^{-1}=m$. Thus the number of self-inverses in 
$GF(p_i^{k_i})$ is the number of elements such that $m\cdot (m-1)\cdot (m+1) = 0$.
Since a field has no zero divisors, these are precisely the elements $0,1$ and $-1$.
Thus if $p_i=2$ and $k_i=1$ then $GF(p_i^{k_i})$ has 2 self-inverses and otherwise 3.\\
(3) Since in a field every element is invertible except $0$, $GF(p_i^{k_i})$ has 
$p_i^{k_i}- 1$ invertibles. $\qed$

\mbox{}\\
{\bf Acknowledgement}: We are indebted to one of the referees of an earlier version of this paper
who observed that our results on invertibles and self-inverses are simple corollaries of the Chinese Remainder Theorem
and ring decomposition.

\end{document}